\documentclass[a4paper]{article}
\usepackage{amsmath,amsthm}

\DeclareMathOperator\id{id}
\DeclareMathOperator\Lin{Lin}

\DeclareMathOperator\SL{SL}
\DeclareMathOperator\SO{SO}
\DeclareMathOperator\SP{Sp}
\DeclareMathOperator\SU{SU}
\DeclareMathOperator\U{U}

\newtheorem{theorem}{Theorem}
\newtheorem{lemma}[theorem]{Lemma}

\newcommand\de{{\rm d}}
\newcommand\e{\varepsilon}
\newcommand\eop{\varepsilon{\circ}\partial}
\newcommand\p{\partial}
\newcommand\x{\times}
\newcommand\ox{\otimes}

\newcommand\dC{{\mathbf C}}

\newcommand\cA{{\cal A}}
\newcommand\cH{{\cal H}}
\newcommand\cO{{\cal O}}
\newcommand\cR{{\cal R}}
\newcommand\cV{{\cal V}}
\newcommand\cX{{\cal X}}

\begin{document}

\title{Construction of Covariant Differential Calculi on 
Quantum Homogeneous Spaces}
\author{Ulrich Hermisson \\
Fachbereich Mathematik, Universit\"at Leipzig \\
Augustusplatz 10, 04109 Leipzig, Germany \\
e-mail: uhermiss@rz.uni-leipzig.de}
\date{October 2, 1998}

\maketitle

{\small {\bf Abstract.} A method of constructing covariant differential
calculi on a quantum homogeneous space is devised. The function algebra
$\cX$ of the quantum homogeneous space is assumed to be a left coideal of a
coquasitriangular Hopf algebra $\cH$ and to contain the coefficients of any
matrix over $\cH$ which is the two-sided inverse of one with entries in
$\cX$. The method is based on partial derivatives. For the quantum sphere of
Podle\'s and the quantizations of symmetric spaces due to Noumi, Dijkhuizen
and Sugitani the construction produces the subcalculi of the standard
bicovariant calculus on the quantum group.

}

\section{Introduction}

The generalization of differential geometric concepts for Lie groups to
quantum groups is ten years after the initiating work of S.~L.~Woronowicz
(cf.~\cite{Wo}) still a promising task. In this paper we propose an
extension of Woronowicz's theory of covariant differential calculi on
quantum groups to a certain class of quantum spaces. The Letter is organized
as follows: Section \ref{genfacts} contains an account of the adopted
notions and notations, a Woronowicz type classification theorem in a dual
formulation and, as guiding example, a new approach to the 3-dimensional
covariant differential calculi on the quantum 2-sphere of Podle\'s
(cf.~\cite{AS,Po1,Po2}). In Section \ref{mainres} the main result, Theorem
\ref{constr}, is proved and quantum spaces of Noumi, Dijkhuizen and Sugitani
(cf.~\cite{Di,NS}) are described as examples. The covariant differential
calculi constructed in these cases are shown to be the subcalculi of the
standard bicovariant calculus on the corresponding quantum group (the chain
rule is valid).

\pagebreak
\section{Definitions and general facts}
\label{genfacts}

We denote by $\cH$ an Hopf algebra over $\dC$ (the complex numbers) with
comultiplication $\Delta$, counit $\e$, antipode $S$ and use Sweedler's
notation: $a_{(1)}\ox a_{(2)} = \Delta(a)$, $a_{(1)}\ox\cdots\ox a_{(Q+1)} =
a_{(1)}\ox\cdots\ox a_{(Q-1)}\ox\Delta(a_{(Q)})$ for $Q = 2,\,3,\,\dots$
Likewise, we use Einstein's convention $a_i\,b^i = \sum_i\,a_i\,b^i$ with
respect to repeated lower case indices.

Let $\cX$ be a subalgebra of $\cH$ and a left coideal, i.~e.
$\Delta(\cX) \subset \cH\ox\cX$, which is what we call a function algebra of
a quantum homogeneous space. Let $(\Gamma,\,\de)$ be a first order
differential calculus on $\cX$, that is, $\Gamma$ an $\cX$-bimodule and
\mbox{$\de:\cX\to\Gamma$} a linear map satisfying 
$\de(a\,b) = \de a\,b + a\,\de b$ (Leibniz rule) such that
\mbox{$\Gamma = \Lin\{a\,\de b\;|\;a,\,b\in\cX\}$.} We require it to be
covariant, i.~e. the linear map $\Delta_\Gamma:\Gamma\to\cH\ox\Gamma: 
a\,\de b\mapsto a_{(1)}\,b_{(1)}\ox a_{(2)}\,\de b_{(2)}$---reflecting the
quantum group operation---to be well-defined. In the following we view
differential calculi only up to isomorphisms.

We fix a basis $e^1,\,\dots,\,e^M$ of a subcomodule of $\cX$. This means
$\Delta(e^i) = \pi^i_j\ox e^j$ with unique $\pi^i_j \in \cH$,
$\Delta(\pi^i_j) = \pi^i_k\ox\pi^k_j$ and $\e(\pi^i_j) = \delta_{ij}$.
In this paper we throughout restrict our considerations to the case
\begin{equation}
\label{freemod}
\de a = a_1\,\de e^1 + \cdots + a_M\,\de e^M
\quad \mbox{ with unique } \quad a_1,\,\dots,\,a_M \in \cX.
\end{equation}
(This assumption serves as a substitute for $\Gamma$ always being a free
left module in the quantum group case.) The formulae
$\de a = \p_i(a)\,\de e^i$ and $\de e^i\,a = \p^{\,i}_j(a)\,\de e^j$ then
define operators on $\cX$ satisfying
$\p^{\,i}_j(a) = \p_j(e^i\,a) - e^i\,\p_j(a)$,
\begin{eqnarray}
\label{prodpd}
& \p_i(a\,b) = \p_j(a)\,\p^{\,j}_i(b) + a\,\p_i(b),
\qquad\p_i(1) = 0, \\[1ex]
& \p^{\,i}_j(a\,b) = \p^{\,i}_k(a)\,\p^{\,k}_j(b),
\qquad\p^{\,i}_j(1) = \delta_{ij}\,1
\end{eqnarray}
as consequences of the Leibniz rule. Thus, the $\p_i$ may be viewed as
generalized skew-derivations which belong to a bialgebra of operators on
$\cX$. If $(\Gamma,\,\de)$ is an inner calculus, that is,
$\de a = \omega\,a-a\,\omega$ for some $\omega = \omega_i\,\de e^i$, then
the formula $\p_i(a) = \omega_j\,\p^{\,j}_i(a) - a\,\omega_i$ presents
$\p_i$ in terms of $\p^{\,j}_i$.

It is important to note that, by covariance, the operators $\p_i$,
$\p^{\,i}_j$ are determined by their values at the quantum group's identity,
which is represented by the counit $\e$ of $\cH$: By covariance
$\Delta_\Gamma(\p_k(a)\,\de e^k)
= (\p_k(a))_{(1)}\,\pi^k_j\ox(\p_k(a))_{(2)}\,\de e^j$
is equal to $\Delta_\Gamma(\de a) = a_{(1)}\ox \de a_{(2)}
= a_{(1)}\ox\p_j(a_{(2)})\,\de e^j$, from this we conclude
$(\p_k(a))_{(1)}\,\pi^k_j\ox(\p_k(a))_{(2)} = a_{(1)}\ox\p_j(a_{(2)})$.
Applying $\id\ox\,\e$, multiplying by $S(\pi^j_i)$ from the right and
dealing with $\p^{\,i}_j$ analogously we get
\begin{eqnarray}
\label{covpd}
& \p_i(a) = a_{(1)}\,\eop_j(a_{(2)})\,S(\pi^j_i), \\[1ex]
\label{covco}
& \p^{\,i}_j(a) = \pi^i_k\,a_{(1)}\,\eop^{\,k}_l(a_{(2)})\,S(\pi^l_j).
\end{eqnarray}
In the case of an inner calculus as before we have
$\eop_i = \e(\omega_j)\,\eop^{\,j}_i - \e(\omega_i)\,\e$.

The restricted dual of an algebra $\cA$, denoted by $\cA^\circ$, by
definition consists of those functionals $f$ on $\cA$ for which finitely
many functionals $f_i,\,f^i$ on $\cA$ exist such that, for all
$a,\,b \in \cA$, the equation $f(a\,b) = f_i(a)\,f^i(b)$ holds. It is
well-known that $\Delta(f) = f_i\ox f^i$ makes $\cA^\circ$ into a coalgebra,
cf.~\cite{KS}.

\begin{theorem}
\label{thclass}
The assignment of the set $\{\eop_1,\,\dots,\,\eop_M\}$ to
$(\Gamma,\,\de)$ establishes a one-to-one correspondence between
\begin{itemize}
\item the covariant first order differential calculi on $\cX$ satisfying
condition \eqref{freemod} and
\item the subsets $\{\chi_1,\,\dots,\,\chi_M\}$ of $\cX^\circ$ for which
$\Lin\{\e,\,\chi_1,\,\dots,\,\chi_M\}$ is a right coideal,
$\chi_i(e^j) = \delta_{ij}$, $\chi_i(1) = 0$ and $\p_i(\cX) \subset \cX$,
where $\p_1,\,\dots,\,\p_M$ are defined by equation \eqref{covpd} with
$\eop_j$ replaced by $\chi_j$.
\end{itemize}
\end{theorem}

\begin{proof}
For a given first order differential calculus on $\cX$ satisfying condition
\eqref{freemod} we apply $\e$ to the first equation of \eqref{prodpd} and
conclude $\eop_i \in \cX^\circ$,
$\Delta(\eop_i) = (\eop_j)\ox(\eop^{\,j}_i) + \e\ox(\eop_i)$. That is,
$\Lin\{\e,\,\eop_1,\,\dots,\,\eop_M\}$ is a right coideal of $\cX^\circ$.
Since $\eop_i(e^j) = \delta_{ij}$ and $\eop_i(1) = 0$, the set
$\{\eop_1,\,\dots,\,\eop_M\}$ has the asserted properties.

Conversely, suppose the set $\{\chi_1,\,\dots,\,\chi_M\}$ complies with the
specified conditions. We define $\p_1,\,\dots,\,\p_M$ as in the statement of
the theorem, thus $\chi_i = \eop_i$. From equation \eqref{covpd} and
condition \eqref{freemod} follows the uniqueness up to isomorphism of the
calculus to which the given set is assigned. To prove the existence, we set
$\Gamma = (\cX\ox\cX) / \Lin\{a\ox b - a\,\p_i(b)\ox e^i\;|\;a,\,b\in\cX\}$
and define by $a\,\de b = p(a\ox b)$ with canonical projection
$p:\cX\ox\cX \to \Gamma$ the left module operation and the map $\de$.
Condition \eqref{freemod} follows from $\de a = \p_i(a)\,\de e^i$ and the
well-definedness of the linear maps $a\,\de b \mapsto a\,\p_i(b)$ on
$\Gamma$, since $\p_i(e^j) = \delta_{ij}\,1$. We make $\Gamma$ into a
bimodule by defining its right module operation through the Leibniz rule
$\de a\,b = \de(a\,b)-a\,\de b$ (note that $\de 1 = 0$ as $\p_i(1) = 0$).
The criterion for this to be well-defined,
$\de a\,b = \p_i(a)\,(\de e^i\,b)$, is equivalent to $\p_i(a\,b) =
\p_j(a)\,\bigl(\p_i(e^j\,b) - e^j\,\p_i(b)\bigr) + a\,\p_i(b)$. This, in
turn, we deduce from $\chi_i(a\,b) = \chi_j(a)\,\bigl(\chi_i(e^j\,b) -
\e(e^j)\,\chi_i(b)\bigr) + \e(a)\,\chi_i(b)$, which ultimately follows from
$\Lin\{\e,\,\chi_1,\,\dots,\,\chi_M\}$ being a right coideal and
$\chi_i(e^j) = \delta_{ij}$, $\chi_i(1) = 0$. Covariance is a consequence of
equation \eqref{covpd}.
\end{proof}

\begin{proof}[Remark]
The right ideal used in \cite{Wo} is the orthogonal complement of the right
coideal used here. Actually any right ideal $\cR$ of $\cX$ with
$\e(\cR) = \{0\}$ gives rise to a covariant first order differential
calculus on $\cX$ as in \cite{Wo}. Indeed, set
$a\,\de b = a\,b_{(1)}\ox p(b_{(2)})$ with $p:\cX\to\cX/\Lin(\cR\cup\{1\})$
the canonical projection, $\Gamma = \Lin\{a\,\de b\;|\;a,\,b\in\cX\}$ and
$\de a\,b = \de(a\,b)-a\,\de b$. It is however an open question whether
this establishes a one-to-one correspondence.
\end{proof}

In the case $\cX=\cH$ we make the more specific assumption
$\de a = a_{ij}\,\de u^i_j$ with unique $a_{ij} \in \cH$,
$i,\,j = 1,\,\dots,\,N$, where $\Delta(u^i_j) = u^i_k\ox u^k_j$ and
$\e(u^i_j) = \delta_{ij}$, and define operators on $\cH$ by
$\de a = \p_{ij}(a)\,\de u^i_j$ and
$\de u^i_j\,a = \p^{\,ij}_{kl}(a)\,\de u^k_l$. Equations \eqref{covpd} and
\eqref{covco} become
\begin{eqnarray}
\label{lcovpd}
& \p_{ij}(a) = a_{(1)}\,\eop_{kj}(a_{(2)})\,S(u^k_i), \\[1ex]
& \p^{\,ij}_{kl}(a) = u^i_m\,a_{(1)}\,\eop^{\,mj}_{nl}(a_{(2)})\,S(u^n_k).
\end{eqnarray}
Bicovariance means that the linear map $a\,\de b\mapsto
a_{(1)}\,b_{(1)}\ox a_{(2)}\,\de b_{(2)} \ox a_{(3)}\,b_{(3)}$ is
well-defined. It can be shown to be equivalent to the condition
\begin{displaymath}
\eop_{kl}(a_{(2)})\,S(a_{(1)})\,a_{(3)} = \eop_{ij}(a)\,S(u^i_k)\,u^l_j
\qquad \mbox{ for all $a$}
\end{displaymath}
in addition to left covariance, and then
$\sigma^{i_1 j_1,\,i_2 j_2}_{k_1 l_1,\,k_2 l_2} =
\eop^{\,i_1 j_1}_{k_2 l_2}(S(u^{i_2}_{k_1})\,u^{l_1}_{j_2})$
commutes with
$\eop^{\,i_1 j_1}_{k_1 l_1}(a_{(1)})\,\eop^{\,i_2 j_2}_{k_2 l_2}(a_{(2)})
\in \dC^{N^4\!\x N^4}$, is invertible and satisfies the braid equation
$(\id\ox\,\sigma)\circ(\sigma\ox\id)\circ(\id\ox\,\sigma)
= (\sigma\ox\id)\circ(\id\ox\,\sigma)\circ(\sigma\ox\id)$ (Woronowicz
braiding).

We will require $\cH$ to be coquasitriangular, i.~e. to be equipped with a
universal $r$-form $r:\cH\ox\cH\to\dC$:
\begin{eqnarray}
\label{cqcc} 
& r(a_{(1)}\ox b_{(1)})\,a_{(2)}\,b_{(2)}
= b_{(1)}\,a_{(1)}\,r(a_{(2)}\ox b_{(2)}), \\[1ex]
\label{rlr}
& r(a\,b\ox c) = r(a\ox c_{(1)})\,r(b\ox c_{(2)}),
\qquad r(1\ox c) = \e(c), \\[1ex]
\label{rrr}
& r(a\ox b\,c) = r(a_{(1)}\ox c)\,r(a_{(2)}\ox b),
\qquad r(a\ox 1) = \e(a).
\end{eqnarray}
The function $\bar r = r\circ(S\ox\id)$ is two-sided convolution inverse
to $r$ (that is, $\bar r(a_{(1)}\ox b_{(1)})\,r(a_{(2)}\ox b_{(2)}) =
\e(a)\,\e(b) = r(a_{(1)}\ox b_{(1)})\,\bar r(a_{(2)}\ox b_{(2)})$). One can
show using \eqref{cqcc} and a left convolution inverse $\bar r$ of $r$ that
$S$ is two-sided composition invertible,
$S^{-1}(a) = r(a_{(1)}\ox S(a_{(2)}))\,S(a_{(3)})
\,\bar r(S^2(a_{(4)})\ox S(a_{(5)}))$, cf.~\cite{KS}.

In the cases $\cH = \cO(\SL_q(N)),\,\cO(\SO_q(N)),\,\cO(\SP_q(N))$,
$q \in \dC\setminus\{0\}$, we denote by $u^i_j$ the canonical generators of
$\cH$. Then $r$ is defined by $r(u^i_k\ox u^j_l) = c\,R^{ij}_{kl}$ if $R$ is
the $R$-matrix of $\cH$, $c = q^{-1/N}$ for $\SL_q(N)$ and $c=1$ else,
cf.~\cite{RTF}.

\subsubsection*{Example 1: Quantum 2-sphere}

Covariant differential calculi on the quantum sphere of Podle\'s which
fulfil condition \eqref{freemod} appear in \cite{AS}. There, in fact, left
and right is reversed, to which however the preceding explanations can be
readily adapted. We view $\cX = \cX_c$ as a right coideal of
$\cH = \cO(\SL_q(2))$, generated as an algebra by $e_{-1}$, $e_0$, $e_1$
with $\Delta(e_i) = e_j\ox\pi^j_i$ and
$c = \e(e_{-1})\,\e(e_1) : \e(e_0)^2$. We define the calculi by requiring
$\de a = \de e_i\,a^i$ with unique $a^i \in \cX_c$ and
\begin{displaymath}
a^i = \p^{\,i}(a) = S(\pi^i_j)\,(\chi^{-1})^{jk}\,\chi_k(a_{(1)})\,a_{(2)},
\end{displaymath}
where the functionals $\chi_k \in \cX_c^\circ$ are specified below and
$\chi^{-1}$ is the inverse of the matrix $\chi$ with coefficients
$\chi_{ij} = \chi_i(e_j)$.

For the parameter value $c = c(3) = -q^6/(q^6+1)^2$ an irreducible
3-dimensional representation $\tau$ of $\cX_c$ exists
(cf.~\cite{Po1} Prop.~4):
\begin{eqnarray*}
\tau(e_{-1}) &=& \e(e_0)\,\frac{q^4-1}{q^6+1}
\,\Biggl(\begin{array}{ccc}
0 & 0 & 0 \\ -1 & 0 & 0 \\ 0 & q^2/(q^2+1) & 0
\end{array}\Biggr), \\
\tau(e_0) &=& \e(e_0)\,\frac{q^4-1}{q^6+1}
\,\Biggl(\begin{array}{ccc}
q^2 & 0 & 0 \\ 0 & q^2-1 & 0 \\ 0 & 0 & -1
\end{array}\Biggr), \\
\tau(e_1) &=& \e(e_0)\,\frac{q^4-1}{q^6+1}
\,\Biggl(\begin{array}{ccc}
0 & -q^2/(q^2+1) & 0 \\ 0 & 0 & q^2 \\ 0 & 0 & 0
\end{array}\Biggr),
\end{eqnarray*}
that is, $\tau_{ij}(e_k) = -\e(e_0)\,\frac{q^8-1}{q^6+1}\,B_{j,ik}$, where
$B_{j,ik}$ is specified in \cite{Po3}. The exceptional calculus for
$c = c(3)$ then emerges from $\chi_i = \e(e_j)\,\tau_{ij} - \e(e_i)\,\e$.

The standard calculus for parameter values $c \neq -q^2/(q^2+1)^2$ results
for $c \neq 0$ from $\chi_i(a) = r(e_i\ox a) - \e(e_i)\,\e(a)$, using the
$r$-form of $\cO(\SL_q(2))$. Theorem~\ref{constr} in the following section
essentially is the statement that this formula and its variant for
$c=\infty$ below generally define such calculi as those considered here. The
matrix $\chi$ is singular exactly if $c \in \{0,\,-q^2/(q^2+1)^2\}$ or
$q^4=1$. The standard calculus and the exceptional one for $c = c(3)$ both
satisfy $\e(e_{-1})\,\e(e_1)\,\de a = -q^4/(q^4-1)^2\,(\omega\,a-a\,\omega)$
with $\omega = (q^2+1)\,\de e_{-1}\,e_1 + \de e_0\,e_0 +
(q^{-2}+1)\,\de e_1\,e_{-1}$.

Finally, we set $\chi_i(a) = r(e_i\ox \nu(a)) - \e(e_i)\,\e(a)$ for the
exceptional calculus in the case $c=\infty$, where $\nu$ is the covariant
algebra automorphism of $\cX_\infty$ with $\nu(e_i) = -e_i$. Here we have
$\de a =
q^4/(\e(e_{-1})\,\e(e_1)\,(q^4+1)\,(q^2+1)^2)\,(\omega\,a-a\,\omega)$.

\subsubsection*{Example 2: Quantum groups}

Let $\cH$ be one of the quantum groups
$\cO(\SL_q(N)),\,\cO(\SO_q(N)),\,\cO(\SP_q(N))$ and set
$\phi^{mn}(a) = r(u^m_l\ox\nu(a)_{(1)})\,r(\nu(a)_{(2)}\ox u^l_n) -
\delta_{mn}\,\e(a)$, using the $r$-form of $\cH$ and a bicovariant algebra
automorphism $\nu$ with $\nu(u^i_j) = \zeta\,u^i_j$, $\zeta\in\dC$,
$\zeta^N=1$ for $\SL_q(N)$ and $\zeta^2=1$ else. We assume that the matrix
$\phi$ with coefficients $\phi^{mn,ij} = \phi^{mn}(u^i_j)$ is invertible;
this excludes a finite set of values for $q$ (e.~g.
$(\zeta^{-1} q^{2/N}\!-1)\,q^2\,\sum_{i=1}^N q^{-2i} \neq q^2-1 \neq 0$ for
$\SL_q(N)$). To obtain the standard bicovariant differential calculi we set
$\eop_{kj} = \phi^{mn}\,(\phi^{-1})_{kj,mn}$ in equation \eqref{lcovpd},
cf.~\cite{HS}.

\section{Construction}
\label{mainres}

We will require $\cX$ to have the following property (for $i$, $j$, $k$ in
any finite set):
\begin{equation}
\label{prop}
\mbox{If }\quad x_{ij} \in \cX,\ y_{ij} \in \cH,
\quad x_{ik}\,y_{kj} = \delta_{ij}\,1 = y_{ik}\,x_{kj},\quad
\mbox{ then }\quad y_{ij} \in \cX.
\end{equation}

\begin{lemma}
\label{lemprop}
Suppose one of the following definitions is applicable to $\cX$:
\begin{enumerate}
\item $\cX = \{a\in\cH\;|\;a_{(1)}\ox p(a_{(2)}) = a\ox p(1)\}$,
where $p$ is a left or right $\cH$-module homomorphism, or
\item $\cX = \{a\in\cH\;|\;a_{(1)}\,\xi(a_{(2)}) = a\,\xi(1)
\mbox{ for all }\xi\in\cV\}$,
where $\cV$ is a linear subspace of $\cH^\circ$ and
$\cV \subset \Delta^{-1}(\cH^\circ\ox\cV)
+ \Delta^{-1}(\cV\ox\cH^\circ)$.
\end{enumerate}
Then $\cX$ has the property \eqref{prop}.
\end{lemma}

\begin{proof}
We suppose that the second case with $\Delta(\cV) \subset \cH^\circ\ox\cV$
is applicable and show that $y_{ij} \in \cX$ for $y_{ij}$ as in
\eqref{prop}. The other cases can be treated similarly. For $\xi\in\cV$ we
have
\begin{eqnarray*}
\lefteqn{y_{ik\,(1)}\,x_{kl\,(1)}\,\xi(y_{ik\,(2)}\,x_{kl\,(2)}) =
y_{ik\,(1)}\,\xi_{(1)}(y_{ik\,(2)})
\,x_{kl\,(1)}\,\xi_{(2)}(x_{kl\,(2)}) =} \hspace{10em} \\
&& y_{ik\,(1)}\,\xi_{(1)}(y_{ik\,(2)})\,x_{kl}\,\xi_{(2)}(1) =
y_{ik\,(1)}\,\xi(y_{ik\,(2)})\,x_{kl}.
\end{eqnarray*}
On the other hand $y_{ik\,(1)}\,x_{kl\,(1)}\,\xi(y_{ik\,(2)}\,x_{kl\,(2)})
= y_{ik}\,x_{kl}\,\xi(1)$, since $y_{ik}\,x_{kl} = \delta_{il}\,1$. We
multiply both results by $y_{lj}$ from the right and use
$x_{kl}\,y_{lj} = \delta_{kj}\,1$ to obtain
$y_{ij\,(1)}\,\xi(y_{ij\,(2)}) = y_{ij}\,\xi(1)$, thus $y_{ij} \in \cX$.
\end{proof}

\begin{proof}[Remark]
That $\cX$ is associated with a quantum subgroup is to say that the first
case of Lemma \ref{lemprop} is applicable with $p$ an Hopf algebra
epimorphism. If $\cH$ is faithfully flat as a left $\cX$-module and $S$ is
bijective, then also the first case is applicable, $p$ being the canonical
projection $\cH \to \cH/\Lin(\cH\,\cX^+)$ with $\cX^+ = (\id-1\,\e)(\cX)$,
cf.~\cite{MW} Thm.~2.1, see also \cite{Br}.
\end{proof}

\begin{theorem}
\label{constr}
Let $\cH$ be a coquasitriangular Hopf algebra, $\cX$ a subalgebra and left
coideal of $\cH$ with property \eqref{prop} and $e^1,\,\dots,\,e^M \in \cX$
with $\Delta(e^i) = \pi^i_j\ox e^j$.
Let $b^1,\,\dots,\,b^M \in \cX$ with $\Delta(b^i) = \psi^i_j\ox b^j$ and
$\nu$ an algebra endomorphism of $\cX$ which is covariant, i.~e.
$\Delta\circ\nu = (\id\ox\,\nu)\circ\Delta$.
Let $\chi^i(a) = r(b^{\,i}\ox\nu(a)) - \e(b^{\,i})\,\e(a)$ and suppose the
matrix $\chi$ with coefficients $\chi^{ij} = \chi^i(e^j)$ is invertible.
Then there uniquely exists a covariant first order differential calculus on
$\cX$ which fulfils condition \eqref{freemod} such that
\begin{displaymath}
\de a = \p_i(a)\,\de e^i, \qquad
\p_i(a) = a_{(1)}\,\chi^k(a_{(2)})\,(\chi^{-1})_{jk}\,S(\pi^j_i).
\end{displaymath}
It satisfies $\de a = \omega\,a-a\,\omega$
with $\omega = \omega_i\,\de e^i$,
$\omega_i = \e(b^k)\,(\chi^{-1})_{jk}\,S(\pi^j_i)$.
\end{theorem}

\begin{proof}
From the first equation in \eqref{rrr} it follows that
$\Lin\{\e,\,\chi^1,\,\dots,\,\chi^M\}$ is a right coideal of $\cX^\circ$,
and $\eop_i(e^j) = \delta_{ij}$, $\eop_i(1) = 0$ hold by construction and
the second equation in \eqref{rrr}. We show $\p_i(\cX) \subset \cX$, in view
of Theorem \ref{thclass} thereby the first assertion is proved. The
equations $a_{(1)}\ox\nu(a_{(2)}) = \nu(a)_{(1)}\ox\nu(a)_{(2)}$, i.~e.
covariance of $\nu$, and \eqref{cqcc} give
\begin{eqnarray*}
\lefteqn{a_{(1)}\,\chi^k(a_{(2)})\,\psi^j_k =
a_{(1)}\,\psi^j_k\,r(b^k\ox\nu(a_{(2)})) - a\,b^j =} \hspace{2em} \\
&& \nu(a)_{(1)}\,\psi^j_k\,r(b^k\ox\nu(a)_{(2)}) - a\,b^j =
r(\psi^j_k\ox\nu(a)_{(1)})\,b^k\,\nu(a)_{(2)} - a\,b^j \in \cX.
\end{eqnarray*}
For $a=e^i$ we get $x^{ij} = \pi^i_l\,\chi^{kl}\,\psi^j_k \in \cX$.
Because of property \eqref{prop} the coefficients of the two-sided inverse
$y_{ji} = S^{-1}(\psi^n_j)\,(\chi^{-1})_{mn}\,S(\pi^m_i)$ are elements of
$\cX$ as well. Thus
$\p_i(a) = a_{(1)}\,\chi^k(a_{(2)})\,\psi^j_k\;y_{ji} \in \cX$. To prove the
second assertion, we note that $\omega_i = b^j\,y_{ji} \in \cX$ and verify
$\p_i(a) = \omega_j\,\p^{\,j}_i(a) - a\,\omega_i$ using the formula
$\p^{\,j}_i(a) = \pi^j_l\,a_{(1)}\,\chi^{nl}\,r(\psi^m_n\ox\nu(a_{(2)}))
\,(\chi^{-1})_{km}\,S(\pi^k_i)$, which can be calculated from
$\p^{\,j}_i(a) = \p_i(e^j\,a) - e^j\,\p_i(a)$.
\end{proof}

\begin{proof}[Remark]
The matrix $\chi$ in the preceding Theorem is surely singular if the
functionals $\chi^1,\,\dots,\,\chi^M$ are not linearly independent, as in
the case of the quantum spheres of Vaksman and Soibelman. This case is
treated in \cite{Schm,We}.
\end{proof}

It is instructive to set up the construction in an alternative way,
resembling that of Jur\v co for quantum groups (cf.~\cite{J}): We retain the
notation of Theorem \ref{constr} and its proof. Let $\Gamma$ be the free
left $\cX$-module with basis $\gamma_1,\,\dots,\,\gamma_M$. We make it into
a bimodule with the right module operation
$\gamma_i\,a = r(\psi^j_i\ox a_{(1)})\,\nu(a_{(2)})\,\gamma_j$. The linear
map $\Delta_\Gamma:\Gamma\to\cH\ox\Gamma:a\,\gamma_i \mapsto
a_{(1)}\,S^{-1}(\psi^j_i)\ox a_{(2)}\,\gamma_j$ then satisfies
$\Delta_\Gamma(\gamma_i\,a) =
S^{-1}(\psi^j_i)\,a_{(1)}\ox\gamma_j\,a_{(2)}$,
and the element $\omega = b^i\,\gamma_i$ is invariant,
i.~e. $\Delta_\Gamma(\omega) = 1\ox\omega$. One readily shows that
$\de a = \omega\,a - a\,\omega$ defines a covariant first order differential
calculus on $\cX$. From the equations $\de e^i = x^{ij}\,\gamma_j$ and
$\gamma_i = y_{ij}\,\de e^j$ one infers that the condition~\eqref{freemod}
is satisfied. The analogous reworking is possible also for the exceptional
calculus for $c = c(3)$ in Example 1 of Section \ref{genfacts}; this leads
to a free right $\cX$-module with basis $\gamma^1,\,\gamma^2,\,\gamma^3$,
left module operation $a\,\gamma^j = \tau_{ij}(a)\,\gamma^i$ and
quantum group operation $\gamma^i\,a \mapsto
\gamma^j\,a_{(1)}\otimes S^{-1}(\pi^i_j)\,a_{(2)}$ (note that
$\pi^i_l\,B_{l,mn} = B_{i,jk}\,\pi^j_m\pi^k_n$).

\subsubsection*{Example 1: Quantum 2-sphere}

The function algebra of the quantum sphere of Podle\'s admits the following
characterization due to Dijkhuizen and Koornwinder (cf.~\cite{DK}):
\begin{displaymath}
\cX_c = \{a \in \cO(\SL_q(2))\;|\;\xi(a_{(1)})\,a_{(2)} = \xi(1)\,a
\mbox{ for all }\xi\in\cV\},
\end{displaymath}
where $\cV = \Lin\{K_{ij}\;|\;i,\,j = 1,\,2\}$,
$K_{ij}(a) = r(a_{(1)}\ox u^l_j)\,J_{kl}\,\bar r(u^k_i\ox a_{(2)})$
with $J = \bigl(\!\begin{smallmatrix}
\e(e_1) \!&\! \e(e_0) \\ 0 \!&\! -q\,\e(e_{-1})
\end{smallmatrix}\!\bigr)$.
This description of $\cX_c$ reduces to the actual condition in~\cite{DK},
``$a\in\cX_c$ iff $\xi(a_{(1)})\,a_{(2)} = 0$'' for a single so-called
twisted primitive element $\xi$. The matrices $J$ are solutions to the
reflection equation with $v^i_j=u^i_j$ specified in Example 2 below,
providing those not taken into account there. Actually the algebra $\cX_c$
properly contains the one generated by $e_{-1}$, $e_0$ and $e_1$ exactly if
$c = -q^n/(q^n+1)^2$, $n$ odd, cf.~\cite{Mue, MS}; for $n=1$ the function
algebra of the quantum hyperboloid (quantum disk) arises.

Since $\Delta(\cV) \subset \cV\ox\cH^\circ$, it follows from a counterpart
of Lemma \ref{lemprop} that $\cX_c$ has the property \eqref{prop}. The
apparent analogue of Theorem \ref{constr} yields the standard 3-dimensional
calculus, if $c \neq 0$, as described in Example 1 of Section
\ref{genfacts}. This is the subcalculus of a bicovariant calculus
(cf.~\cite{Po2} Thm.~1(a) and \cite{AS}, but the argument given in the
following Example, suitably adapted, applies as well).

\subsubsection*{Example 2: Quantum symmetric spaces}

Noumi, Dijkhuizen and Sugitani (cf.~\cite{Di,NS}) have defined
quantizations of the classical irreducible compact Riemannian symmetric
spaces
\begin{eqnarray*}
&& \SU(N)/\SO(N),\quad\SO(N)/\U(N/2),\quad\SO(N)/(\SO(L)\x\SO(N-L)), \\
&& \SU(N)/\SP(N),\quad\SP(N)/\U(N/2),\quad\SP(N)/(\SP(L)\x\SP(N-L)), \\
&& (*)\;\U(N)/(\U(L)\x\U(N-L)).
\end{eqnarray*}
They assign to each of these spaces an invertible solution
$J \in \dC^{N\x N}$ to the reflection equation
\begin{displaymath}
r(u^j_a\ox u^i_b)\,J^{ac}\,r(v^k_c\ox u^b_d)\,J^{dl} =
J^{ia}\,r(v^c_a\ox u^j_b)\,J^{bd}\,r(v^l_c\ox v^k_d),
\end{displaymath}
$v^i_j=u^i_j$ or in case $(*)$ $v^i_j=S(u^j_i)$, and define the function
algebra of the quantum space by
\begin{displaymath}
\cX = \{a\in\cH\;|\;a_{(1)}\,\xi(a_{(2)}) = a\,\xi(1)
\mbox{ for all }\xi\in\cV\},
\end{displaymath}
where $\cV = \Lin\{K^{ij}\;|\;i,\,j = 1,\,\dots,\,N\}$,
$K^{ij}(a) = \bar r(a_{(1)}\ox u^i_k)\,J^{kl}\,r(v^j_l\ox a_{(2)})$.
Since $\Delta(\cV) \subset \cV\ox\cH^\circ$, Lemma \ref{lemprop} implies
that $\cX$ has the property \eqref{prop}.

Let $e^1,\,\dots,\,e^M$ be a basis of the subcomodule
$(\id-1\,h)\bigl(\Lin\{u^i_k\,J^{kl}\,v^j_l\;|
\;i,\,j = 1,\,\dots,\,N\}\bigr)$ of $\cX$. The Haar functional $h$ on $\cH$
used here is uniquely determined by $h(1) = 1$ and one of the equations
$a_{(1)}\,h(a_{(2)}) = h(a)\,1 = h(a_{(1)})\,a_{(2)}$ (cf.~\cite{KS}; its
application here is to split off the $1$, if necessary). For some of the
spaces it has been shown that these elements generate $\cX$ as an algebra
(cf.~\cite{DN} Prop.~3.11 and \cite{No} Thm.~4.3).

We prove that if Theorem \ref{constr} with $b^i=e^i$ and $\nu=\id$ yields a
calculus on $\cX$ it is the subcalculus of the bicovariant calculus with
$\nu=\id$ specified in Example 2 of Section \ref{genfacts}: For $a\in\cX$
we have
\begin{eqnarray*}
\phi^{mn}(a)
&=& r(u^m_l\ox a_{(1)})\,r(a_{(2)}\ox u^l_i)
\,K^{ij}(1)\,(J^{-1})_{jn} - \delta_{mn}\,\e(a) \\
&=& r(u^m_l\ox a_{(1)})\,r(a_{(2)}\ox u^l_i)
\,K^{ij}(a_{(3)})\,(J^{-1})_{jn} - \delta_{mn}\,\e(a) \\
&=& r(u^m_l\,J^{li}\,v^j_i\ox a)\,(J^{-1})_{jn} - \delta_{mn}\,\e(a),
\end{eqnarray*}
hence $\phi^{mn}(a) = c^{mn}_k\,\chi^k(a)$ with
$c^{mn}_k = \phi^{mn}(e^i)\,(\chi^{-1})_{ik}$. From this we conclude
$\p_{ij}(a) = \p_k(a)\,\p_{ij}(e^k)$. The calculation also shows
$\chi^k(a) = d^k_{mn}\,\phi^{mn}(a)$ with $d^k_{mn} \in \dC$, from which we
deduce that $a_k\,\p_{ij}(e^k) = 0$ implies $a_k=0$. This completes the
proof.

Moreover, by the above calculation, instead of proving that the matrix
$\chi$ in Theorem \ref{constr} is nonsingular one may equivalently verify
that the rank of the matrix with coefficients
$c^{mn,ij} = \phi^{mn}(u^i_k\,J^{kl}\,v^j_l)$ is $M$. The result is that
$\chi$ is nonsingular exactly if
$(q^{4/N}\!-1)\,q^2\,\sum_{i=1}^N q^{-2i} \neq q^4-1 \neq 0$ for
$\SU(N)/\SO(N)$, where $M = N\,(N+1)/2$, and
$(q^{4/N}\!-1)\,q^4\,\sum_{i=1}^{N/2} q^{-4i} \neq q^2-1 \neq 0$ for
$\SU(N)/\SP(N)$, where $M = N\,(N-1)/2$. The remaining cases still have to
be examined.

\subsubsection*{Acknowledgements}

I am very indebted to Prof.~K.~Schm\"udgen for his encouragement to
investigate covariant differential calculi on quantum spaces. I thank
I.~Heckenberger for stimulating discussions. This work was supported by the
Deutsche Forschungsgemeinschaft within the scope of the postgraduate
scholarship programme ``Graduiertenkolleg Quantenfeldtheorie'' at the
University of Leipzig.

\end{document}